\documentclass[12pt]{amsart}
\usepackage{amsmath,amssymb,amsfonts,amsthm}
\usepackage{enumerate}
\usepackage{hyperref}

\title[Characterizing the powerset]{Characterizing the powerset by a complete (Scott) sentence}

\author[I. Souldatos]{Ioannis  Souldatos}
\address{4001 W.McNichols, Mathematics Department, University of Detroit Mercy, Detroit, MI 48221, USA}
\email{souldaio@udmercy.edu}


\date{}

\newcommand{\omegaone}{\ensuremath{\omega_1}}
\newcommand{\lomegaone}{\ensuremath{\mathcal{L}_{\omega_1,\omega}}}
\newcommand{\lomega}{\ensuremath{\mathcal{L}_{\omega,\omega}}}
\newcommand{\alephs}[1]{\ensuremath{\aleph_{#1}}}
\newcommand{\alephalpha}{\alephs{\alpha}}
\newcommand{\alephomegaone}{\alephs{\omegaone}}
\newcommand{\alephalphaplus}{\alephs{\alpha+1}}
\newcommand{\beths}[1]{\ensuremath{\beth_{#1}}}

\newcommand{\bethomegaone}{\beths{\omegaone}}

\newcommand{\M}{\ensuremath{\mathcal{M}}}
\newcommand{\N}{\ensuremath{\mathcal{N}}}
\newcommand{\A}{\ensuremath{\mathcal{A}}}
\newcommand{\B}{\ensuremath{\mathcal{B}}}
\newcommand{\F}{\ensuremath{\mathcal{F}}}

\newcommand{\C}{\ensuremath{\mathcal{C}}}

\newcommand{\G}{\ensuremath{\mathcal{G}}}

\newcommand{\lang}[1]{\ensuremath{\mathcal{L}_{#1}}}

\newcommand{\ch}{\ensuremath{\mathcal{CH}_{\omega_1,\omega}}}
\newcommand{\homch}{\ensuremath{\mathcal{HCH}_{\omega_1,\omega}}}
\newcommand{\ltok}{\ensuremath{\lambda^{\kappa}}}
\newcommand{\ltoomega}{\ensuremath{\lambda^{\omega}}}

\newcommand{\f}[1]{\ensuremath{f(#1)}}
\newtheorem{theorem}{Theorem}[section]
\newtheorem{lemma}[theorem]{Lemma}
\newtheorem{corollary}[theorem]{Corollary}
\theoremstyle{definition}
\newtheorem{definition}[theorem]{Definition}
\newtheorem{claim}[theorem]{Claim}
\newtheorem{fact}{Fact}
\newtheorem*{conjecture}{Conjecture}

\newtheorem{observation}[theorem]{Observation}
\newtheorem{property}{Property}

\newtheorem{remark}[theorem]{Remark}

\numberwithin{equation}{section}

\begin{document}


\baselineskip=17pt


\begin{abstract} This paper is part II of a study on cardinals that are characterizable by a Scott sentence, continuing the work in \cite{CharacterizableCardinals}. A cardinal $\kappa$ is characterized by a
Scott sentence $\phi_\M$, if $\phi_\M$ has a model of size $\kappa$, but no model of size $\kappa^+$.

The main question in this paper is the following: Are the characterizable cardinals closed under
 the powerset operation? We prove that if $\alephs{\beta}$ is characterized by a Scott sentence, then $2^{\alephs{\beta+\beta_1}}$ is (homogeneously) characterized by a Scott sentence, for all $0<\beta_1<\omegaone$ (cf. theorem  \ref{mainhompowersetthrm} and definition \ref{defhch}). So, the answer to the above question is positive, except the case $\beta_1=0$ which remains open. 

As a consequence we derive that if $\alpha\le\beta$ and $\alephs{\beta}$ is characterized by a Scott sentence, then $\alephs{\alpha+\alpha_1}^{\alephs{\beta+\beta_1}}$ is  (homogeneously) characterized by a Scott sentence, for all  $\alpha_1<\omegaone$ and $0<\beta_1<\omegaone$ (cf. theorem \ref{miscpowers}).
Whence, depending on the model of ZFC, we see that the class of characterizable and homogeneously characterizable cardinals
 is much richer than previously known. 
Several open questions are mentioned at the end.
\end{abstract}

\maketitle


\section{Introduction}
This paper is part II of a study on cardinals that are characterizable by a Scott sentence. We refer the reader to \cite{CharacterizableCardinals} for more details and background information. The main question we try to answer in this paper is the following: Are the characterizable cardinals closed under  the powerset operation? 

We prove a positive answer for all cardinals of the form $2^{\aleph_{\beta+\beta_1}}$, where $\aleph_{\beta}$ is characterized by a Scott sentence and $0<\beta_1<\omegaone$. The case $\beta_1=0$ remains open. The main construction is contained in theorem \ref{powersetofsuccessor}: Given a cardinal $\aleph_{\beta}$ that is characterized by a Scott sentence, we prove that $2^{\aleph_{\beta+1}}$ is also characterized by a Scott sentence. The idea is to create a complete graph whose edges are $\aleph_{\beta+1}$-colored and which \lq\lq{}mimics\rq\rq{} the behavior of $2^{\aleph_{\beta+1}}$ (cf. property \ref{crucialproperty}). This will ensure that the graph can not have size greater than $2^{\aleph_{\beta+1}}$ and the theorem follows. 

Throughout the whole paper we work only with countable languages, but we do not restrict ourselves to relational languages only.

\subsection*{Basic Definitions}
We start by mentioning the basic definitions. 

\begin{definition}\label{charcarddef} We say that a \lomegaone-sentence $\phi$ \emph{characterizes}
\alephalpha, or that $\alephalpha$ is\emph{ characterizable}, if $\phi$ has models in all cardinalities up to $\alephalpha$, but not in cardinality
\alephalphaplus. If $\phi$ is the Scott sentence of a
countable model (or any other complete sentence), we say that it
\emph{completely characterizes} \alephalpha, or that \alephalpha\; is
\emph{completely characterizable}. Denote by \ch, the set of all
completely characterizable cardinals.
\end{definition}

W. Hanf proved in \cite{HanfsIncompactness} that there exists a cardinal,
denoted by $\mathcal{H}(\lomegaone)$, such that if an $\lomegaone$- sentence has a
model in  cardinality $\mathcal{H}(\lomegaone)$, then it has models in \emph{all}
cardinalities. $\mathcal{H}(\lomegaone)$ is called
the Hanf number for $\lomegaone$ and it is equal to
\bethomegaone. So, $\ch\subset\bethomegaone$ and from now on we only consider cardinals below
\bethomegaone. We will also restrict ourselves to cardinals that are completely characterizable, and we may refer to them as just characterizable cardinals.

\begin{definition}\label{defhch} If $P$ is a unary predicate symbol, we say that it is
\emph{completely homogeneous} for the \lang{}- structure \A, if
$P^{\A}=\{a|\A\models P(a)\}$ is infinite and every permutation of
it extends to an automorphism of \A.

If $\kappa$ is a cardinal, we will say that $\kappa$ is
\emph{homogeneously characterizable} by
$(\phi_{\kappa},P_{\kappa})$, if $\phi_{\kappa}$ is a complete
\lomegaone- sentence and  $P_{\kappa}$ a unary predicate in the
language of $\phi_{\kappa}$ such that
\begin{itemize}
    \item $\phi_{\kappa}$ doesn't have models of power $>\kappa$,
    \item if \M\; is the (unique) countable model of
    $\phi_{\kappa}$, then $P_{\kappa}$ is infinite and completely
    homogeneous for \M\; and
    \item there is a model \A\; of $\phi_{\kappa}$ such that
    $P_{\kappa}^{\A}$ has cardinality $\kappa$.
\end{itemize}

Denote the  set of all homogeneously characterizable cardinals by \homch. 
\end{definition}
Obviously, $\homch\subset\ch$ and in \cite{HjorthsKnightPaper} Hjorth proved that $\aleph_0\in\ch\setminus\homch$. So, $\homch$ is a strict subset of $\ch$ and it is open whether there is another example in $\ch\setminus\homch$ besides $\aleph_0$. 

In the same paper \cite{HjorthsKnightPaper} Hjorth also proved
that the class of characterizable cardinals is closed under
successors and countable unions, i.e. if $\alephalpha\in\ch$ and
$\beta<\omegaone$, then $\alephs{\alpha+\beta}\in\ch$. This means
that characterizable cardinals come into clusters of length
$\omegaone$.
\begin{definition}\label{headofcluster} A cardinal $\alephs{\alpha}\in\ch$ is called the \emph{head of a
cluster}, if we can not find ordinals $\beta,\gamma$ such that
\begin{itemize}
    \item $\alephs{\gamma}\in\ch$,
    \item $\beta<\omegaone$ and
    \item $\alephalpha=\alephs{\gamma+\beta}$
\end{itemize}
\end{definition}
It is immediate that all characterizable cardinals are of the form
$\alephs{\alpha+\beta}$, where $\alephs{\alpha}$ is the head of a
cluster and $\beta<\omegaone$.

\section{Fraisse Construction}\label{FraisseConstruction}

We describe briefly a Fraisse-type construction which we will use in section \ref{powerset}. Definition \ref{fingenoverM} is an expanded version of a similar definition in \cite{LinearOrderings} and the theorems that follow definition \ref{fingenoverM} are from \cite{LinearOrderings} too. The interested reader should refer there for more details on the proofs.\footnote{As mentioned in \cite{LinearOrderings} too, these theorems are straightforward generalizations of Fraisse\rq{}s theorems in the context of ``finitely generated over $\M$ structures''. Hjorth in \cite{HjorthsKnightPaper} used similar constructions as well as Laskowski and Shelah in \cite{LaskowskiShelahsOnTheExistenceOfAtomicModels}. The proofs of theorems  \ref{Fraisse} and \ref{Fraisse2} follow the proofs of the original Fraisse theorems. So, our effort was in organizing into a unique framework tools that other authors have used before.}.

In definition \ref{fingenoverM} we work with a relational vocabulary $\tau$. If $\tau$ is not relational and $\M$ is a $\tau$- structure, replace all constant symbols $c\in\tau$ by corresponding (unary) predicate  symbols $P_c$ and all function symbols $f$ of arity $n$ by predicate symbols $P_f$ of arity $n+1$ and stipulate: 
\begin{enumerate} 
\item $P_c^\M(x)$ iff $c^\M=x$ and
\item $P_f^\M(\vec{x},y)$ iff $f^\M(\vec{x})=y$. 
\end{enumerate}
So, without loss of generality assume that $\tau$ is a relational vocabulary.

\begin{definition}\label{fingenoverM} Fix a relational vocabulary $\tau$ and a $\tau$-model $\M$. Let  $M$ be a unary predicate not in $\tau$ and $\star$ be a constant symbol not in $\tau$. Abusing notation, $\star$ denotes both the constant symbol and its interpretation, which we now fix to be some element not in $\M$.
\begin{enumerate}
\item Let $\tau'$  be a (not necessarily relational) vocabulary such that $\tau'\supset\tau\cup\{M,\star\}$ and  let $\A\supset(\M\cup\{\star\})$ be a $\tau'$- model. Call $\A$ a model over $\M$ (write ``model/$\M$'') if all the following hold:
\begin{enumerate}
\item for any $x\in\A$, $\A\models M(x)$ iff $x\in\M$. 
\item if $R\in\tau$, then the interpretation of $R$ is retained on $\M$ and is void outside of $\M$. I.e. for any tuple $\vec{x}\in\A$, $\A\models R(\vec{x})$ iff $\vec{x}\in \M$ and $\M\models R(\vec{x})$, and
\item if $f$ is a function symbol in $\tau'\setminus \tau$, then $f$ is of interest only when defined on tuples $\vec{a}\in \A\setminus\M$. If $\vec{a}$ contains an element in $\M$, then let $f(\vec{a})=\star$. Other than that, we make no other use of the constant symbol $\star$.
\end{enumerate} 

\item Let $\A$ be a model/$\M$ and $A_0\subset \A$. Denote by $<A_0>$ the $\tau'$- substructure of $\A$ generated by $A_0$. If $B=<A_0\cup M>$, where $A_0$ is a finite subset of $\A\setminus M$, then $B$ will be called \emph{finitely generated over} $\M$. Abbreviate  \lq\lq{}finitely generated over $\M$\rq\rq{} as \lq\lq{}finitely generated/$\M$\rq\rq{}. Similar abbreviations apply below too. 

\item If $\A_0\subset \A_1$ are two models/$\M$, and $\A_0$ is a $\tau'$- substructure of $\A_1$, then $\A_0$ will be called a \emph{substructure over $\M$} of $\A_1$, written as $\A_0\subset_\M \A_1$. If $f:\A_0\rightarrow\A_1$,  $f|_\M=id_\M$ and $f(\A_0)\subset_\M \A_1$, then $f$ will be called an \emph{embedding} over $\M$ of $\A_0$ into $\A_1$. Again, we have  \lq\lq{}substructure/$\M$\rq\rq{} and \lq\lq{}embedding/$\M$\rq\rq{} being the obvious abbreviations. 

\item If $\A_0,\A_1$ are two models/$\M$, then $\A_0,\A_1$ are \emph{isomorphic over} $\M$, denoted by $\A_0\cong_\M \A_1$, if there exists an isomorphism $i:\A_0\rightarrow \A_1$ such that $i|_{\M}=id_{\M}$. If $\A_0=\A_1$, then every isomorphism over $\M$ will be called an \emph{automorphism} over $\M$. A  partial function $i$ between $\A_0$ and $\A_1$ is a  \emph{partial isomorphism} over $\M$, if $\M\subset dom(i)$ and $i$ is an isomorphism/$\M$ between $dom(i)$ and $range(i)$. \lq\lq{}Isomorphism/$\M$\rq\rq{}, \lq\lq{}automorphism/$\M$\rq\rq{} and ``partial isomorphism/$\M$'' are again the obvious abbreviations.

\item If $\A_0,\A_1$ are two models/$\M$, they will be called ``back-and-forth-equivalent over $\M$'' (abbreviate as ``back and forth equivalent/$\M$'') if there is a back-and-forth system of partial isomorphisms/$\M$ between $\A_0$ and $\A_1$.

\item If in the usual Joint Embedding Property (JEP) and Amalgamation Property (AP) ``embeddings'' are replaced by ``embeddings/$\M$'', we get the Joint Embedding Property over $\M$ (JEP/$\M$) and the  Amalgamation Property over $\M$ (AP/$\M$) respectively. 
\end{enumerate}
\end{definition}

\textbf{Conventions and Notation:}
\begin{enumerate}
\item For the rest of the paper, JEP and AP denote JEP/$\M$ and AP/$\M$ respectively.
\item Let $M|_\tau\models\phi$ denote the fact that for all $\tau'$ structures $\A$ the interpretation of the predicate $M$ in $\A$ restricted to the vocabulary $\tau$ satisfies the sentence $\phi$. If $\phi$ is an $\lomegaone$ sentence, then $M|_\tau\models\phi$ can be expressed by an $\lomegaone$ sentence. If $\A$ is a $\tau'$ sturcture, then for the rest of the paper $M(\A)=\M$ and $M(\A)\models\phi$ will denote that $M(\A)|_\tau=\M$ and $M(\A)|_\tau\models\phi$ respectively, although some times we will include the subscripts for emphasis.
\end{enumerate}

The key observation is that Fraisse's theorems hold even in the context of finitely generated/$\M$ substructures.

\begin{theorem}\label{Fraisse} Fix a countable model $\M$ and let $K(\M)$ be a countable collection of finitely generated/$\M$ structures (up to isomorphism/$\M$). If $K(\M)$ has the Hereditary Property (HP), the Joint Embedding Property (JEP) and the  Amalgamation Property (AP), then there is a countable structure $\F$ which we will call the \emph{Fraisse limit} of $K(\M)$ such that
\begin{enumerate}
\item $\F$ is a model/$\M$ and is unique up to isomorphism/$\M$ ,
\item $K(\M)$ is the collection of all finitely generated/$\M$ substructures of $\F$ (up to isomorphism/$\M$), and
\item every isomorphism/$\M$ between finitely generated/$\M$ substructures of $\F$ extends to an automorphism/$\M$ of $\F$.
\end{enumerate}

The converse is also true, i.e. if $\F$ is a countable model/$\M$ such that  every isomorphism/$\M$ between finitely generated/$\M$ substructures of $\F$ extends to an automorphism/$\M$ of $\F$ and $K(\M)$ is the collection of all finitely generated/$\M$ substructures of $\F$, then $K(\M)$ has the HP, the JEP and the AP.

We will use the notation $\lim K(\M)$ when referring to the Fraisse limit.
\begin{proof}  Both theorems \ref{Fraisse} and \ref{Fraisse2} follow from the corresponding proofs of the original theorems of Fraisse by making the obvious modifications. One can consult  \cite{HodgesShorterModelTheory} for more details.
\end{proof}
\end{theorem}

\begin{theorem}\label{Fraisse2} Fix a model $\M$. Assume that $\A,\B$ are two (not necessarily countable) models/$\M$  such that
\begin{itemize}
\item for every finitely generated/$\M$ substructures $C\subset D$ of $\A$ (or of $\B$), and every embedding/$\M$ $f:C\mapsto \A$ ($f:C\mapsto \B$), there is an embedding/$\M$ $g:D\mapsto \A$ ($g:D\mapsto \B$) that extends $f$, and
\item the collection of all finitely generated/$\M$ substructures of $\A$ (up to isomorphism/$\M$) is the same as the collection of  all finitely generated/$\M$ substructures of $\B$ (up to isomorphism/$\M$).
\end{itemize}
Then $\A$ and $\B$ are back-and-forth equivalent/$\M$. In particular, $\A$ and $\B$ are back-and-forth equivalent and $\A\equiv_{\infty,\omega}\B$.
\end{theorem}

We now give a slightly different version of theorem \ref{Fraisse} that will be more fitting to work with in the next section.

\begin{theorem} \label{propertiesIandII} Fix a countable model $\M$ and let $\phi_\M$ be its Scott sentence. Let $\psi$ be an $\lomegaone(\tau')$- sentence all whose models are finitely generated substructures over some model of $\phi_\M$\footnote{This can be expressed by $M|_\tau\models\phi_\M$.}. Let $K(\M)$ be the collection (up to isomorphism/$\M$) of all finitely generated/$\M$ structures that satisfy $\psi$. If $K(\M)$ is countable and satisfies the HP, the JEP and the AP, then the Fraisse limit $\F$ of $K(\M)$ satisfies the conjunction of
\begin{description}
  \item[(I)$_{\M}$] Every finitely generated/$\M$ substructure satisfies $\psi$.
  \item[(II)$_{\M}$] For every  $A_0$ finitely generated/$\M$ substructure, if $A_1\supset A_0$ and $A_1$ satisfies $\psi$, then there exists some finitely generated/$\M$ substructure $B$ and an isomorphism $i:B\cong A_1$, such that $A_0\subset B$ and $i|_{A_0}=id$.
\end{description}
Moreover, the conjunction of (I)$_{\M}$ and (II)$_{\M}$ can be written as an \lomegaone-sentence and if $\A$ is a countable model/$\M$ that satisfies this conjunction, then $\A$ is isomorphic/$\M$ to $\F$.
\begin{proof} The fact that $\F$ satisfies the conjunction of (I)$_{\M}$ and (II)$_{\M}$ follows from theorem \ref{Fraisse}. The proof that  the conjunction $(I)_{\M}\wedge (II)_{\M}$ can be written as an \lomegaone-sentence is straightforward, given the assumptions on $\psi$. The proof that $\A$ is isomorphic to $\F$ is a standard back-and-forth (over $\M$) argument.
\end{proof}
\end{theorem}

\begin{corollary} \label{uniquelimit} Fix a countable model $\M$, let $\phi_\M$ be its Scott sentence and let $\M'$ be another countable model of $\phi_\M$. Let $\psi$ be an $\lomegaone(\tau')$- sentence all whose models are finitely generated substructures over some model of $\phi_\M$. For any $\N$ model of $\phi_\M$, let $K(\N)$ be the collection of all finitely generated/$\N$ structures (up to isomorphism/$\N$) that satisfy $\psi$.

Assume that $K(\M)$ is countable and satisfies the HP, the JEP and the AP. Then:
\begin{enumerate}[(a)]
\item $\lim K(\M)\cong \lim K(\M')$ and
\item if $\A$ is a countable model/$\M'$ that satisfies $(I)_{\M'}\wedge (II)_{\M'}$, then $\A\cong\lim K(\M)$.
\end{enumerate}
\begin{proof} (a) It is immediate by the assumptions that $K(\M')$ is also countable and satisfies the HP, the JEP and the AP.

Let $i$ witness that $\M\cong\M'$ and let $\tilde{i}$ be the following mapping:

$\tilde{i}(x)=
\begin{cases}
x & \text{, if $x\in \lim K(\M)\setminus\M$}\\
i(x)& \text{, if $x\in \M$}
\end{cases}
$

Obviously, $\tilde{i}(\lim K(\M))$ is a countable structure/$\M'$ that satisfies both $(I)_{\M'}$ and $(II)_{\M'}$. By theorem \ref{propertiesIandII}, it is isomorphic to $\lim K(\M')$ and part (a) follows. 

(b) Immediate from part (a) and theorem \ref{propertiesIandII}.
\end{proof}
\end{corollary}

We need a bit more. Assume $\M,\phi_\M,\M',\psi,K(\M)$ are as in corollary \ref{uniquelimit} and consider the following variation of $(I)$ and $(II)$, the difference than before being that finitely generated substructures are over the (interpretation of the) predicate $M$.
\begin{description}
  \item[(I)] Every finitely generated/$M$ substructure satisfies $\psi$.
  \item[(II)] For every  $A_0$ finitely generated/$M$ substructure, if $A_1\supset A_0$ and $A_1$ satisfies $\psi$, then there exists some finitely generated/$M$ substructure $B$ and an isomorphism $i:B\cong A_1$, such that $A_0\subset B$ and $i|_{A_0}=id$.
\end{description}
It is immediate that any model of $(I)_{\M}\wedge (II)_{\M}$ will also satisfy $(I)\wedge (II)$ and the same is true for all models of $(I)_{\M'}\wedge (II)_{\M'}$. By corollary \ref{uniquelimit}(b), the conjunction of $(I)\wedge (II)$ together with $M|_\tau\models\phi_\M$ yields a complete sentence. In addition, it is expressible in $\lomegaone$ that a substructure is finitely generated/$M$, which makes $(I)\wedge (II)\wedge M|_\tau\models\phi_\M$ an $\lomegaone$- sentence. Consequently, all its models (countable or uncountable) are $\equiv_{\omegaone,\omega}$- equivalent to each other.

\begin{corollary}\label{EquivalentToLimKM} Let $\M,\phi_\M,\psi,K(\M)$ be as in corollary \ref{uniquelimit} and let $\N$ be another model of $\phi_\M$ (possibly uncountable). If $\A$ is a model/$\N$ that satisfies $(I)\wedge (II)$, then $\A$ is  $\equiv_{\infty,\omega}$- equivalent to $\lim K(\M)$.
\begin{proof}  By the comments before the corollary. 
\end{proof}
\end{corollary}

If $\M$ is a countable model whose Scott sentence $\phi_\M$ characterizes a certain cardinal $\kappa$, we will use the Scott sentence of $\lim K(\M)$ to characterize some cardinal $\lambda\ge\kappa$. In order to construct (uncountable) models of the Scott sentence of $\lim K(\M)$ we will use the following theorem.

\begin{theorem} \label{towardsfullness}  Let $\M,\phi_\M,\psi,K(\M)$ be as in corollary \ref{uniquelimit} and assume that $\phi_\M$ characterizes an infinite cardinal $\kappa$, $\N$ is a model of $\phi_\M$ and $\lambda$ is a cardinal $\ge\kappa$. Moreover, assume that:
\begin{enumerate}
  \item The number of the finitely generated substructures/$\N$ (up to isomorphism/$\N$) that satisfy $\psi$ is $\le\lambda$. In particular, if $A$ is a finitely generated/$\N$ structure, then there are $\le\lambda$ many (non-isomorphic/$\N$) structures/$\N$ that satisfy $\psi$ and extend $A$, and
  \item For every $\G$ a model/$\N$  such that
     \[|\G\setminus M(\G)|<\lambda,\; \G\mbox{ satisfies }(I)_{\N}\]

and for any $A_0,A_1$ finitely generated/\N\; structures with \[A_0\subset\G,\; A_1\supset A_0, \mbox{ and }A_1 \text{ satisfies $\psi$},\]

there exists some model/$\N$ $\G'\supsetneqq \G$ such that 
\[ |\G'\setminus M(\G)|<\lambda,\; \G'\mbox{ satisfies }(I)_{\N}\]

and there is some finitely generated/\N\; structure $B\subset\G'$ and an isomorphism $i:B\cong A_1$, with $A_0\subset B$ and $i|_{A_0}=id$.
\end{enumerate}
Under the assumptions 1 and 2, we conclude that there exists a model/$\N$ $\G^*$ with $|\G^*\setminus M(\G)|=\lambda$ and $\G^*$ satisfies $(I)_{\N}$ and $(II)_{\N}$. By theorem \ref{EquivalentToLimKM}, $\G^*$ also satisfies the Scott sentence of $\lim K(\M)$.
\end{theorem}
Notice that in part (1) of theorem \ref{towardsfullness}, the number of non-isomorphic structures in $K(\N)$ is allowed to be strictly less than $\lambda$. Actually, in the construction of section \ref{powerset} the number of non-isomorphic structures in $K(\N)$ will be $\kappa^+$, while the size of $\G^*$ will be $2^{(\kappa^+)}$, i.e. $\lambda=2^{(\kappa^+)}$.  
Also, notice that the structure $\G\rq{}$ in part (2) is required to be \emph{strictly} bigger than $\G$.

\section{Known Theorems}\label{KnownTheorems}
 
This section contains certain known theorems about characterizable cardinals. They are quoted from \cite{BaumgartnersHanfNumber}, \cite{HjorthsKnightPaper}, \cite{CharacterizableCardinals} and \cite{LinearOrderings}.

\begin{theorem}[theorem 3.4 from \cite{CharacterizableCardinals}]\label{ltoomegahom} If $\lambda\in\ch$, then $\ltoomega\in\homch$.
\end{theorem}

\begin{corollary}[corollary 3.6 from \cite{CharacterizableCardinals}]\label{ltokhom} If $\kappa$ is an infinite cardinal and $\ltok\in\ch$, then $\ltok\in\homch$.
\end{corollary}

\begin{theorem}[theorem 3.7 from \cite{CharacterizableCardinals}]\label{headpower} If $\alephalpha^{\alephs{\beta}}\in\ch$, then for all $\gamma<\omegaone$,
\[\alephs{\alpha+\gamma}^{\alephs{\beta}}\in\homch.\]
\end{theorem}

From  \cite{BaumgartnersHanfNumber} we have the following theorem, which also appears in \cite{CharacterizableCardinals}
\begin{theorem}[Baumgartner]\label{baumgartner}
If $\kappa\in\homch$, then $2^{\kappa}\in\homch$.
\end{theorem}

The following two theorems can be derived from \cite{HjorthsKnightPaper} and  appear in \cite{CharacterizableCardinals} and  \cite{LinearOrderings}\footnote{The formulation of theorem \ref{loofsizekappa+} as it appears here is slightly different than the one in \cite{LinearOrderings}.} respectively.

\begin{theorem}[Hjorth]\label{limitthrm} Whenever $\alephs{\alpha_n}$, $n\in\omega$, is an non-decreasing sequence of cardinals in \ch,  then $\alephs{\lambda}= \sup{\alephs{\alpha_n}}$ is also in \ch.
\end{theorem}

\begin{theorem}(Hjorth) \label{loofsizekappa+} If $\kappa\in\ch$, then at least one of the following holds:
\begin{enumerate}
  \item $\kappa^+\in\homch$ or,
  \item there is a countable model $\M$ in a language that contains a unary predicate $P$ and a binary predicate $<$ and whose Scott sentence $\phi_0$
  \begin{enumerate}
    \item has no models of cardinality $\kappa^{++}$, 
    \item does have a model of cardinality $\kappa^+$,
    \item in every model of $\phi_0$, $<$ is a dense linear order without endpoints and
    \item in every model $\N$ of $\phi_0$ of size $\kappa^+$,  every initial segment of $(P^\N,<^\N)$ has size $\le\kappa$.
  \end{enumerate}
  \end{enumerate}
\end{theorem}

With all these in place we are ready to prove some new theorems.

\section{Powersets}\label{powerset}

\begin{theorem}\label{powersetofsuccessor} If $\kappa\in\ch$, then $2^{(\kappa^+)}\in\homch$.
\end{theorem}
We will split the proof into a series of lemmas and theorems.  By corollary \ref{ltokhom}, if $2^{(\kappa^+)}\in\ch$, then $2^{(\kappa^+)}$ is also in $\homch$. So, it suffices to prove $2^{(\kappa^+)}\in\ch$. 

By theorem \ref{loofsizekappa+} we have to consider two cases. If $\kappa^+\in\homch$ the result is immediate from theorem \ref{baumgartner}. If this is not the case, then there exists a countable model $\M$ as in part (2) of theorem \ref{loofsizekappa+}. I.e. if $\phi_0$ is the Scott sentence of $\M$, then the language of $\phi_0$ contains a unary predicate $P$ and a binary predicate $<$ and if $\N$ is any model of $\phi_0$, then
\begin{enumerate}
  \item[(a)] $|\N|\le\kappa^+$,
  \item[(b)] $(P^\N,<^\N)$ is a dense linear ordering without endpoints that has size $\le\kappa^+$ and
  \item[(c)] for every $m\in P^\N$, the initial segment $\{m_1\in P^\N|m_1<^\N m\}$ has size $\le\kappa$.
\end{enumerate}
Moreover, there is a certain model $\N$ for which equality holds true in all the above properties. 

\begin{definition} \label{NotationOfN} For the rest of the proof, fix $\M$ and $\phi_0$ as in part (2) of theorem \ref{loofsizekappa+}, and fix $\N$ that makes equality true for (a),(b) and (c) above. Let  $\tau$ be the signature of $\phi_0$ and it follows that $\tau\supset\{P,<\}$. 
\end{definition}

The following is immediate from the definition of $\N$ and it will be used later in the proof.

\begin{fact}\label{CofinalFact} An increasing sequence in $P(\N)$ is cofinal iff it has length $\kappa^+$. 
\end{fact}

Our goal is to construct a certain $\psi\in\lomegaone$ and let $K(\M)$ be the collection of all finitely generated/$\M$ structures that satisfy $\psi$. By lemma \ref{JEPandAP}, $K(\M)$ will satisfy both the JEP and the AP, while the HP is immediate. By theorem \ref{propertiesIandII} $K(\M)$ has a (Fraisse) limit, which we will call $\F$. The Scott sentence of $\F$, called $\phi_\F$, witnesses that $2^{(\kappa^+)}\in\ch$, i.e. it does not have any models of size $>2^{(\kappa^+)}$ (theorem \ref{NotTooBig}), but it does have a model of size equal to $2^{(\kappa^+)}$ (theorem \ref{SomeBig}). This last assertion is established by considering $K(\N)$, the collection of all finitely generated structures/$\N$ that satisfy $\psi$, and using theorem \ref{towardsfullness}.

The idea behind the construction is to try to mimic the behavior of the powerset. If $\kappa$ is an infinite cardinal and $a,b\in 2^\kappa$, then let $f(a,b)$ be equal to the least $\alpha\in\kappa$ such that $a(\alpha)\neq b(\alpha)$. It is immediate that all distinct $a,b,c\in 2^\kappa$ satisfy the following property which we isolate.

\begin{property}\label{crucialproperty} For distinct $a,b,c\in 2^\kappa$, exactly two of $f(a,b),f(a,c)$ and $f(b,c)$ are equal, while the third one is larger than the other two.\end{property}

Property \ref{crucialproperty} can be expressed in $\lomegaone$ (indeed in $\lomega$) and is the one that drives the whole construction. We will consider all elements of $2^\kappa$ as vertices of a complete graph $V$ and we will color all the edges between them using colors from a linearly ordered set $M$ in such a way that property \ref{crucialproperty} is satisfied. The symbol for the coloring function will be $f$.

\begin{definition} 
Let $\tau$ be the signature of $\phi_0$ fixed by definition \ref{NotationOfN}. Without loss of generality we assume that $\tau$ is a relational signature. If this is not the case, work similarly to the comments before definition \ref{fingenoverM}. Let $V,M,f,\star$ be new symbols not in $\tau$ such that $V(\cdot), M(\cdot)$ are new unary predicate symbols, $f$ is a new binary function symbol and $\star$ is a new constant symbol. Let $\tau\rq{}=\tau\cup\{V,M,f\}$.
\end{definition}
We will consider $\tau\rq{}$- structures  $\A$  such that:
\begin{itemize}
\item $\A\supset\M$ ($\M$ fixed by definition \ref{NotationOfN}) and $\A$ is a ``model over $\M$'' (cf. definition \ref{fingenoverM}). In particular, $M(\A)=\M$, $M(\A)\models \phi_0$ and there are predicate symbols $P$ and $<$ in $\tau$ such that $P(\A)\subset M(\A)$ and $<^\A$ is a linear order on $P(\A)$ which is void outside $M(\A)$. 
\item $V(\A)$ and $M(\A)$ (the interpretations of $V$ and $M$ in $\A$ respectively) partition the whole universe.
\end{itemize}
We make the following convention.

\textbf{Convention}: $P(\A)=M(\A)$. In general, $P(\A)$ is a subset of $M(\A)$. Under this assumption, the whole $M(\A)$ becomes linearly ordered by $<^\A$.\footnote{The convention that $P(\A)=M(\A)$ is made for convenience. There is no problem if we assume that $P(\A)$ is a strict subset of $M(\A)$, but in this case we would have to do a little more work; namely in clause (2) of definition \ref{PsiDefinition}, $f(x,y)$ must be an element of $P$ (not $M$) for all $x,y$. In order to simplify the notation and do not carry both predicates $M$ and $P$ around.}

\begin{definition}\label{PsiDefinition} Let $\psi$ be the conjunction of the following: 

\begin{enumerate}
  \item $V\cup M$ partition the universe, $V$ is finite and $M|_\tau\models\phi_0$. In particular, there is a linear order $<$ defined on a subset of $M$ and following the above convention, we will assume that $<$ is a linear order on the whole $M$. $V$ will be seen as a set of vertices in a complete graph, while $M$ will be an ordered set of colors. Colors will be assigned to edges.
  \item For all $x,y$, $f(x,y)\in M$ and $f(x,y)=f(y,x)$. If $\Delta=\{(v,v)|v\in V\}$, then $f$ is of interest only on domain $V^2\setminus\Delta$ and assume that for all $(x,y)\notin V^2\setminus\Delta$, $f(x,y)=\star$, where $\star$ is a fixed value (cf.  definition \ref{fingenoverM} too). For $a\neq b\in V$, we think of $f$ as assigning to the edge $(a,b)$ a color in $M$.
  \item For every distinct $a_0,a_1,a_2\in V$, if $\f{a_0,a_1}\neq \f{a_0,a_2}$, then
  \[\f{a_1,a_2}=\min \{\f{a_0,a_1},\f{a_0,a_2}\}.\]
  Otherwise, $\f{a_1,a_2}> \f{a_0,a_1}=\f{a_0,a_2}$. 
\end{enumerate}

Let $K(\M)$ be the collection of all $\tau'$-structures $\A$ (up to isomorphism/$\M$) such that $\A$ is a model over $\M$ (cf.  definition \ref{fingenoverM}) and $\A\models\psi$. In particular, $M(\A)=\M$. 
\end{definition}
Notice that clause (3) is a reformulation of property \ref{crucialproperty} and it is the one that will require work to prove. Also, since the relations from $\tau$ are void outside of $\M$, the isomorphism type of any structure $\A$ in $K(\M)$ is determined by the values of $f^\A$ on the finite set $V(\A)$. Hence, there are only $|\M|^{<\omega}=|\M|$ many possibilities for $f$ and we conclude:

\begin{observation} \label{NumberOfIsomorhpicCopies} The size of $K(\M)$ is the same as the size of $\M$ and the same holds true if we replace $\M$ by any other model of $\phi_0$ (possibly uncountable).
\end{observation} 

The following is also immediate. 

\begin{observation} \label{ReformulationOfClause3}
If $a,b,c$ are distinct elements in $V$, and $\f{a,b}=m$, then at
least one of $\f{a,c},\f{b,c}$ is $\le m$.
\end{observation}

\begin{lemma}\label{JEPandAP} $K(\M)$ satisfies the J.E.P. and the A.P.
\begin{proof}
For  the A.P., let $\A,\B,\C\in K(\M)$, with $\A\subset\B,\C$. In particular, $M(\A)=M(\B)=M(\C)=\M$. We can assume that $\B\cap\C=\A$ and we aim at making the union of $\B$ and $\C$ to a $\tau'$- structure in $K(\M)$. This entails two steps: (1) Make $\B\cup\C$ into a $\tau'$- structure by defining $f$ on pairs of the form $(b,c)$, $b\in\B$, $c\in\C$, and (2) prove that under this definition $\B\cup\C$ satisfies $\psi$. In particular, prove clause (3) of definition \ref{PsiDefinition}. Without loss of
 generality we can assume that \[|V(\B)\setminus V(\A)|=|V(\C)\setminus V(\A)|=1.\]
If this is not the case, use (double) induction on $|V(\B)\setminus V(\A)|$ and $|V(\C)\setminus V(\A)|$.

So, assume that $V(\B)\setminus V(\A)=\{b\}$ and that
$V(\C)\setminus V(\A)=\{c\}$. We define $\f{b,c}$ by cases: If there exists some $a\in V(\A)$ such that $\f{a,b}\neq\f{a,c}$, then let $\f{b,c}=\min\{\f{a,b},\f{a,c}\}$.
If this is not the case, then for all $a\in V(\A)$ it holds true that $\f{a,b}=\f{a,c}$. Let $m$ be a value in $\M$ which is greater than all elements in $ range(f^\B)\cup range(f^\C)$ and assign $\f{b,c}=m$. 

A couple of comments before we proceed: 
\begin{enumerate}[(a)]
\item In the above definition $\f{a,b}$ means $f^\B(a,b)$ and $\f{a,c}$ means $f^\C(a,c)$. Since $f^\A$, $f^\B$ and $f^\C$ agree on all elements of $V(\A)$, superscripts do not make any difference and will be omitted. 
\item Since $V(\B)$ and $V(\C)$ are finite, the same is true for $range(f^\B)$ and $range(f^\C)$ and the existence of $m$ in the second case is guaranteed by the fact that $(\M,<)$ is a dense linear order without endpoints. 
\end{enumerate}
We now prove that
\begin{claim}\label{claim1} The value of $\f{b,c}$ is well-defined.
\begin{proof}
Assume that there exists some $a_1\in V(\A)$ such that $\f{a_1,b}\neq\f{a_1,c}$ and let $m_1=\min\{\f{a_1,b},\f{a_1,c}\}$. We have to prove that for all other values $a_2\in V(\A)$, if $\f{a_2,b}\neq\f{a_2,c}$, then \[\min\{\f{a_2,b},\f{a_2,c}\}=m_1.\]

Towards contradiction assume that there is some $a_2\in V(\A)$ such that $\f{a_2,b}\neq\f{a_2,c}$ and $\min\{\f{a_2,b},\f{a_2,c}\}\neq m_1$. Let $m_2$ be $\min\{\f{a_2,b},\f{a_2,c}\}$. We will derive a contradiction in the case that $m_2> m_1$ and  $m_1=\f{a_1,b}$. The proof in all other cases is done similarly. 

Thus, $a_1,a_2,b$ are all elements in $V(\B)$ and it holds that $\f{a_1,b}=m_1<m_2\le \f{a_2,b}$. By assumption, $\B\models\psi$, which implies that \[\f{a_1,a_2}=
\min\{\f{a_1,b},\f{a_2,b}\}=\f{a_1,b}=m_1.\] By observation \ref{ReformulationOfClause3}, $\f{a_1,a_2}=m_1$  implies that one of $\f{a_1,c},\f{a_2,c}$ has value $\le m_1$. But
$m_1<m_2\le \f{a_2,c}$. So, it must be that $\f{a_1,c}\le m_1$. If $\f{a_1,c}=m_1$, then $\f{a_1,b}=\f{a_1,c}=m_1$ which contradicts the assumption that $\f{a_1,b}\neq \f{a_1,c}$.
If $\f{a_1,c}<m_1$, then $\f{a_1,c}<\min\{\f{a_1,b},\f{a_1,c}\}$, which is also a contradiction. So, in either case we get a contradiction, which proves that $\f{b,c}$ is well- defined.
\end{proof}
\end{claim}

\begin{claim}\label{claim2} Under the above definition of $\f{b,c}$, $\B\cup\C\models\psi$. In particular, $\B\cup\C\in K(\M)$. 
\begin{proof}
Clauses (1) and (2) of definition \ref{PsiDefinition} are immediate and only clause (3) requires some work. I.e. we have to prove that for all $a_1\in V(\A)$, if $\f{a_1,b}\neq\f{a_1,c}$, then $\f{b,c}=\min\{\f{a_1,b},\f{a_1,c}\}$. Otherwise, $\f{b,c}>\f{a_1,b}=\f{a_1,c}$.
We split the proof into two cases:

Case I: There is some $a\in V(\A)$ such that $\f{a,b}\neq \f{a,c}$. Then $\f{b,c}=\min\{\f{a,b},\f{a,c}\}$. Let $a_1\neq a$ be a different element in $V(\A)$. If $\f{a_1,b}\neq \f{a_1,c}$, then by the proof of the previous claim, $\f{b,c}=\min \{\f{a_1,b},\f{a_1,c}\}$, which gives the result. If $\f{a_1,b}=\f{a_1,c}$, we will prove that $\f{b,c}>\f{a_1,b}=\f{a_1,c}$. Towards contradiction assume that $\f{b,c}\le\f{a_1,b}=\f{a_1,c}$. Without loss of generality also assume that $\min\{\f{a,b},\f{a,c}\}=\f{a,b}$ and the other case is proved similarly. We split the proof into two sub-cases:

Subcase I.1: $\f{b,c}=\f{a_1,b}=\f{a_1,c}$. It follows that $\f{a,b}=\min\{\f{a,b},\f{a,c}\}=\f{b,c}=\f{a_1,b}=\f{a_1,c}$. Since $\B\models\psi$ and $a,a_1,b$ are all elements in $\B$ such that $\f{a,b}=\f{a_1,b}$, we conclude that $\f{a,a_1}>\f{a,b}=\f{a_1,b}$. Moreover, $\f{a,c}>\f{a,b}=\min\{\f{a,b},\f{a,c}\}=\f{b,c}=\f{a_1,c}$. Since $\C\models\psi$ and $a,a_1,c$ are all elements in $\C$ such that $\f{a,c}>\f{a_1,c}$, we conclude that $\f{a,a_1}=\min\{\f{a,c},\f{a_1,c}\}=\f{a_1,c}$. Combining these two results together, $\f{a,a_1}>\f{a_1,b}=\f{a_1,c}=\f{a,a_1}$. Contradiction.

Subcase I.2: $\f{b,c}<\f{a_1,b}=\f{a_1,c}$. Again, it follows that $\f{a,b}=\min\{\f{a,b},\f{a,c}\}=\f{b,c}<\f{a_1,b}=\f{a_1,c}$ and since $\B\models\psi$ and $a,a_1,b$ are all elements in $\B$ such that $\f{a,b}<\f{a_1,b}$, we conclude that $\f{a,a_1}=\min\{\f{a,b},\f{a_1,b}\}=\f{a,b}$. On the other hand, since $\C\models\psi$ and $a,a_1,c$ are elements in $\C$, by observation \ref{ReformulationOfClause3}, one of $\f{a,c}$ and $\f{a_1,c}$ must be $\le\f{a,a_1}$. But $\f{a,c}>\min\{\f{a,b},\f{a,c}\}=\f{a,b}=\f{a,a_1}$ and $\f{a_1,c}>\f{b,c}=\f{a,b}=\f{a,a_1}$. So, both cases give a contradiction. 

Case II: For all $a_1\in V(\A)$, $\f{a_1,b}=\f{a_1,c}$.
 Then $\f{b,c}$ is bigger than any value in $range(f^\B)\cup range(f^\C)$. In particular, $\f{b,c}>\f{a_1,b}=\f{a_1,c}$, for all $a_1\in V(\A)$, which concludes the proof.
\end{proof}
\end{claim}

Finally, for the J.E.P., let  $\A=\B\cap\C$ and work as in the proof of the A.P. If $V(\A)=\emptyset$, while $V(\B)=\{b\}$ and $V(\C)=\{c\}$, then $range(f^\B)\cup  range(f^\C)$ is empty and assign $\f{b,c}$ to be \emph{any} value. The details follow.
\end{proof}
\end{lemma}
The fact that $K(\M)$ satisfies the H.P. is immediate and by theorem \ref{Fraisse}, the Fraisse limit of $K(\M)$ exists. 
\begin{definition}\label{FraisseDefinition} Let $\F$ be the Fraisse limit of $K(\M)$ (cf. theorem \ref{Fraisse}), where $K(\M)$ is the collection given by definition \ref{PsiDefinition}. Then $\F$ is countable and unique up to isomorphism/$\M$. 

Also let $\phi_\F$ be the Scott sentence of $\F$. By theorem \ref{propertiesIandII}, $\phi_\F$ is equivalent to  $(I)\wedge(II)\wedge M|_\tau\models\phi_\M$.
\end{definition}
By corollary \ref{uniquelimit}, the choice of $\M$ is not of importance and it can be replaced by any other (isomorphic) model of $\phi_0$. However, if $\G$ is a (possibly uncountable) model of $\phi_\F$, then $M(\G)\models\phi_0$, and if $\G$ is countable, then $M(\G)\cong\M$. The Scott sentence $\phi_\F$ is the one that witnesses that $2^{(\kappa^+)}\in\ch$, i.e. it does not have any models of size $>2^{(\kappa^+)}$ (theorem \ref{NotTooBig}), but it does have a model of size equal to $2^{(\kappa^+)}$ (theorem \ref{SomeBig}).

\begin{theorem}\label{NotTooBig} If $\G\models\phi_\F$, then $|\G|\le 2^{(\kappa^+)}$.
\begin{proof} Towards contradiction assume that $\G\models\phi_\F$ and $|\G|> 2^{(\kappa^+)}$. Without loss of generality, $|\G|= (2^{\kappa^+})^+=(\beths{1}(\kappa^+))^+$. Since $M(\G)\models\phi_0$ and $\phi_0$ characterizes $\kappa^+$, 
 it must be $|M(\G)|\le \kappa^+$, which further implies that
 $|V(\G)|=(\beths{1}(\kappa^+))^+$.
 Hence, $f^\G$ gives a function from $[(\beths{1}(\kappa^+))^+]^2$ to
 $\kappa^+$.

 By Erd\"{o}s-Rado theorem, \[\beths{1}(\kappa^+)^+\rightarrow (\kappa^{++})^2_{\kappa^+},\]
 which also implies
 \[\beths{1}(\kappa^+)^+\rightarrow (3)^2_{\kappa^+}.\]
 Applying this to $f^\G$ we get that there is a homogeneous
 set of size $3$, i.e. there are $a,b,c\in V(\G)$ such that
 \[f^\G({a,b})=f^\G({b,c})=f^\G({a,c}),\] which contradicts the way \F\; was defined and the
 fact that $\G\models\phi_\F$.
\end{proof}
\end{theorem}

\begin{theorem}\label{SomeBig} There exists some $\G^*\models\phi_\F$ such that $|\G^*|=2^{(\kappa^+)}$.
\begin{proof} Let $\N$ be the model of $\phi_0$ fixed by definition \ref{NotationOfN}. In particular, \begin{enumerate}
  \item[(a)] $|\N|=\kappa^+$,
  \item[(b)] $(P^\N,<^\N)$ is a dense linear ordering without endpoints that has size equal to $\kappa^+$ and
  \item[(c)] for every $m\in P^\N$, the initial segment $\{m_1\in P^\N|m_1<^\N m\}$ has size equal to $\kappa$.
\end{enumerate}
Following the convention right before definition \ref{PsiDefinition} we assume that $<^N$ is a linear order on the whole $\N$. 

Let $K(\N)$ be the collection of all $\tau'$-structures $\A$ (up to isomorphism/$\N$) such that $\A$ is a model over $\N$ and $\A\models\psi$ (cf. definition \ref{PsiDefinition}). By observation \ref{NumberOfIsomorhpicCopies} the number of non-isomorphic/$\N$ structures in $K(\N)$ is $\kappa^+$ and we will prove that $K(\N)$ also satisfies assumption (2) of theorem \ref{towardsfullness} for $\lambda=2^{(\kappa^+)}$, which will conclude the proof.

Thus, assume that $\G$ is a model over $\N$ such that

$$|\G|< 2^{(\kappa^+)}\text{ and }\G\text{ satisfies }(I)_{\N}.$$

Also assume that $\A\subset \G$, $\A\subset \B$, $\A,\B\in K(\N)$. Then, we have to prove that there is a structure $\G\rq{}\supsetneqq\G$ that satisfies $(I)_{\N}$, and there is some $\C\subset \G'$ and $i:\B\cong\C$ with $i|_\A=id$.

Without loss of generality assume that $V(\B)\setminus V(\A)=\{b\}$. If this is not the case,  proceed  by induction on $|V(\B)\setminus V(\A)|$. We can also assume that $b\notin\G$. If this is not the case, replace $\B$ by an isomorphic copy. Let $\G\rq{}=\G\cup\{b\}$ and define $f^{\G\rq{}}$ to agree with $f^\G$ on $\G$ and to agree with $f^{\B}$ on $\B$. This is well-defined, since $f^\G$ and $f^\B$ agree on $\A$, and with the exception of a couple of instances, we will drop superscripts from $f$ for the rest of the proof. 

The goal is to:
\begin{enumerate}
\item extend the definition of $f^{\G\rq{}}$ on pairs of the form $(b,v)$ for all $v\in V(\G)\setminus V(\A)$. Then $\G'$ becomes a $\tau'$- structure and
\item prove that any finitely generated/$\N$ substructure of $\G\rq{}$ is in $K(\N)$. Equivalently, any finitely generated/$\N$ substructure of $\G\rq{}$ satisfies $\psi$. Clauses (1) and (2) of $\psi$ are immediate, but clause (3) takes some work again.
\end{enumerate}

Part (1) of the goal is accomplished in two steps by definitions \ref{VomegaDefinition} and \ref{FDefinition}. The latter definition involves certain sequences that we call \emph{selectors} (cf. definition \ref{SelectorDefinition}) and we will prove first some needed lemmas about selectors before we can define $f^{\G\rq{}}(b,v)$, $v\in V(\G)\setminus V(\A)$.

Part (2) of the goal is the content of lemma \ref{PsiForGPrime}.

We begin with the following observation: Let $b$ the unique element of $V(\B)\setminus V(\A)$ and $v$ some element in $V(\G)\setminus V(\A)$. Notice that the value $f(b,v)$ has not been defined yet. If there exists some $a\in V(\A)$ such that $f(a,b)\neq f(a,v)$, then the only way to define $f(b,v)$ so that the substructure/$\N$ generated by $a,b,v$ satisfies sentence $\psi$, is to set $f(b,v)=\min\{f(a,b),f(a,v)\}$. The point is that we do not have freedom to select the value of $f(b,v)$. We can actually push this argument even further applying the same observation $\omega$ many times. 

\begin{definition}\label{VomegaDefinition} Define the following sets:
\begin{itemize}
  \item $V_0=V(\A)$.
  \item $V_{n+1}=\{v\in V(\G)| \exists a\in V_n (\f{a,b}\neq \f{a,v})\}$, for all $n\in\omega$.
\item $V_\omega=\cup_n V_n$.
\item $W_0=V(\G)\setminus V_\omega$
\end{itemize}
If $v\in V_{n+1}$, let $\f{b,v}=\min \{\f{a,v},\f{a,b}\}$, where $a$ is given by the fact that $v\in V_{n+1}$. 
\end{definition}

Under this definition we have extended the definition of $f^{\G'}$ on the whole $V_\omega$. We need to prove that it is well-defined and it actually satisfies $\psi$ as we intended to.

\begin{claim} The value of $\f{b,v}$, for $v\in V_\omega$, is well- defined. 
\begin{proof} By induction on $n$ and working as in the proof of claim \ref{claim1} of lemma \ref{JEPandAP}. The details are left to the reader.
\end{proof}
\end{claim} 

\begin{claim}\label{VOmegaClaim} Any finitely generated/$\N$ substructure of $V_\omega\cup \{b\}$ satisfies sentence $\psi$ (cf. definition \ref{PsiDefinition}). In particular, there is an isomorphic/$\N$ copy of it in $K(\N)$.
\begin{proof} Again, by induction on $n$ and working this time as in the proof of claim \ref{claim2} of lemma \ref{JEPandAP}.
\end{proof}
\end{claim}

It remains to define $\f{b,w}$ for $w\in W_0$.

We first observe that 
\[\tag{*}\label{star}W_0=\{w\in V(\G)| \forall a\in V_\omega,\;  \f{a,b}=\f{a,w}\}.\] 

We work similarly as before: Let $b$ be the unique element of $V(\B)\setminus V(\A)$, $w$ some element in $W_0$ and the value $f(b,w)$ has not been defined yet. Since for all $a\in V_\omega$  $f(a,b)=f(a,w)$, then the only way to define $f(b,w)$ so that the substructure/$\N$ generated by $a,b,w$ satisfies sentence $\psi$, is to set $f(b,w)>f(a,b)=f(a,w)$. I.e. every $a\in V_\omega$ is imposing a
restriction on each $w\in W$ that $\f{w,b}$ must be defined greater than $\f{b,a}=\f{a,w}$.  On
the face of it, we have many restrictions to worry about, but in fact there are only finitely many.

\begin{claim}\label{VZeroClaim} If $m\in M(\G)$ is greater than $\f{b,a}$, for all $a\in V(\A)=V_0$, then it is greater than $\f{b,v}$, for all $v\in V_\omega=\cup_n V_n$.
\begin{proof} By induction on $n$. Assume that this is true for $V_n$ and that $v\in V_{n+1}$ as witnessed by $a\in V_n$, i.e. $\f{a,b}\neq\f{a,v}$. Then, $\f{b,v}=\min\{\f{a,b},\f{a,v}\}\le \f{a,b}<m$.
\end{proof}
\end{claim}

The key point here is that these finitely many restrictions
do not even depend on $w$. 
\begin{definition} Let $m_{-1}$ to be the value \[\max\{\f{b,a}|a\in V_0\}.\] 
\end{definition} 

By \eqref{star} it follows that 
\[ m_{-1}=\max\{\f{a,w}|a\in V_0,\; w\in W_0\}\] 
and we have to ensure that for all $w\in W_0$, $\f{b,w}$ is defined greater than $m_{-1}$. 
Before we define $f(b,w)$ (definition \ref{FDefinition}), we give some motivation first: Let $v_0\in W_0$ be an arbitrary element, let  $m_0$ be some value greater than $m_{-1}$ and define $f(b,v_0)=m_0$. Once $f(b,v_0)$ has been defined, if $w$ is an element such that $f(w,v_0)\neq m_0$, then, arguing as before, the substructure/$\N$ generated by $b,w,v_0$  satisfies $\psi$ only if we define $f(b,w)=\min\{f(b,v_0)=m_0,f(w,v_0)\}$. The only elements for which this argument does not apply are the elements of the set $A_1:=\{v\in W_0|f(v,v_0)=m_0\}$. If $v$ is an element in $A_1$, then $f(b,v_0)=m_0=f(v_0,v)$ and the substructure/$\N$ generated by $b,v_0,v$ satisfies $\psi$ only if $f(b,v)$ is defined greater than $m_0$. We can continue in this fashion replacing $W_0$ by $A_1$ and $m_{-1}$ by $m_0$. Then select some $v_1\in A_1$ and some $m_1>m_0$ and define $f(b,v_1)=m_1$ and $f(b,w)=\min\{f(b,v_1),f(v_1,w)\}$ for all $w\in A_1$ such that $f(v_1,w)\neq m_1$ and carry the elements of $A_2:=\{v\in A_1|f(v,v_1)
=m_1\}$ in the next stage. Repeat this process $\kappa^+$ many times, selecting at every stage some $m_\alpha$ greater than all previous $m_\beta$, $\beta<\alpha$, some $v_\alpha\in A_\alpha$ and letting $A_{\alpha+1}:=\{v\in A_\alpha|f(v,v_\alpha)=m_\alpha\}$. We stop when  we either find some $\alpha$ such that $A_\alpha=\emptyset$, in which case we have defined $f^{\G'}$ on all pairs $(b,v)$, $v\in V(\G)\setminus V(\A)$, and $\G'$ is now a $\tau'$- structure, or we complete $\kappa^+$ many steps, in which case the $m_\alpha$'s form a cofinal sequence in $\N$ (cf. fact \ref{CofinalFact}) and there is no room for repeating the process further. If after $\kappa^+$ many steps the intersection of all $A_\alpha$, $\alpha<\kappa^+$, is empty, then we have defined $f^{\G'}$ on all pairs $(b,v)$, $v\in V(\G)\setminus V(\A)$, making $\G'$ to a $\tau'$- structure again. Otherwise there is at least one element $w\in W_0$ such that $w\in A_\alpha$, for all $\alpha<\kappa^+$, and it is impossible to define $f(b,w)$ to 
be greater than all $m_\alpha$, $\alpha<\kappa^+$. Therefore, this last case must be avoided and the motivation behind selectors is to do exactly that. 

More formally now:

\begin{definition}\label{SelectorDefinition} A sequence $\vec{A}=\langle A_\alpha,v_\alpha,m_\alpha\rangle$ of length $\kappa^+$ will be called a \emph{selector}, if

\begin{itemize}
  \item $m_\alpha\in M(\G)=\N$ and $\beta<\alpha\rightarrow m_\beta<m_\alpha$, 
  \item $A_0=W_0$,
  \item $A_\gamma=\cap_{\alpha<\gamma} A_\alpha$, for $\gamma$ limit ordinal, 
  \item if $A_{\alpha}=\emptyset$, then $A_{\alpha+1}=\emptyset$ and there is no restriction on $v_\alpha$, and 
  \item if $A_{\alpha}\neq\emptyset$, then $v_\alpha$ must be an element of $A_\alpha$ and 
  \[A_{\alpha+1}= \{v\in A_\alpha| v\neq v_\alpha\text{ and }\f{v,v_\alpha}=m_\alpha\}.\]
\end{itemize}
The selector will be called \emph{good}, if
\[\bigcap_{\alpha<\kappa^+} A_\alpha=\emptyset.\]
Otherwise, it will be called \emph{bad}.
\end{definition}
It follows from the definition that the $A_\alpha$\rq{}s form a decreasing sequence of subsets of $V(\G)$, while all the $m_\alpha$'s form an increasing sequence in $M(\G)$. We continue until we get the empty set, or we complete $\kappa^+$ many steps.

The existence of a good selector is what we are after.

\begin{lemma} For all selectors $|\bigcap_{\alpha<\kappa^+} A_\alpha|\le 1$.
\begin{proof} This is obvious for good selectors, so let  $\vec{A}$ be a bad selector with $a_0,a_1\in \bigcap_{\alpha<\kappa^+} A_\alpha\subset V(\G)$ and $\f{a_0,a_1}=m$. Since $(m_\alpha)_{\alpha<\kappa^+}$ is an increasing sequence of length $\kappa^+$, by Fact \ref{CofinalFact}, this sequence is cofinal and we can choose some $\alpha<\kappa^+$ such that $m_\alpha>m$.
Then, $a_0,a_1\in A_{\alpha+1}$, which implies that \[\f{a_0,v_\alpha}=\f{a_1,v_\alpha}=m_\alpha.\] 
Since every finitely generated/$\N$ substructure of $\G$ satisfies $\psi$, this is true in particular for the structure generated by $a_0,a_1,v_\alpha$. Therefore, by clause (3) of $\psi$ it must be that $\f{a_0,a_1}>\f{a_0,v_\alpha}=\f{a_0,v_\alpha}$. This is a contradiction, since $\f{a_0,a_1}=m<m_\alpha=\f{a_0,v_\alpha}=\f{a_1,v_\alpha}$.

Hence, $|\bigcap_{\alpha<\kappa^+} A_\alpha|= 1$, for all bad selectors.
\end{proof}
\end{lemma}

\begin{definition} For a bad selector $\vec{A}$, call $v(\vec{A})$ the unique element in $\bigcap_{\alpha<\kappa^+} A_\alpha$.

For selectors $\vec{A},\vec{B}$ (good or bad), let $\vec{A}\sim\vec{B}$ if and only if
\[\bigcap_{\alpha<\kappa^+} A_\alpha=\bigcap_{\alpha<\kappa^+} B_\alpha.\]
For bad selectors, this is the case iff $v(\vec{A})=v(\vec{B})$. If $\vec{A}$ is a good selector, then $\bigcap_{\alpha<\kappa^+} A_\alpha=\emptyset$ and we conclude that all good selectors are equivalent.
\end{definition}

This defines an equivalence relation on the selectors, each one associated with a (unique) element in $V(\G)$ or the empty set. So, it is immediate that there are exactly $|V(\G)|<2^{\kappa^+}$ many equivalence classes. We will use this fact to prove that there is at least one good selector. We need the following easy lemma before that:

\begin{lemma} Let $\alpha_0<\kappa^+$ and $\vec{A}=\langle A_\alpha,v_\alpha,m_\alpha\rangle$, $\vec{B}=\langle B_\alpha,w_\alpha,m'_\alpha\rangle$ be bad selectors such that:
\begin{itemize}
\item for all $\alpha\le\alpha_0$, $A_\alpha= B_\alpha\neq\emptyset$, 
\item $v_{\alpha_0}=w_{\alpha_0}$ and
\item $m_{\alpha_0}\neq m'_{\alpha_0}$.
\end{itemize}
Then $\vec{A}\nsim\vec{B}$.
\begin{proof} Towards contradiction, assume that $v=v(\vec{A})=v(\vec{B})$. Then we get on the one hand that $\f{v,v_{\alpha_0}}=m_{\alpha_0}$ and on the other hand $\f{v,w_{\alpha_0}}=m'_{\alpha_0}$. Contradiction.
\end{proof}
\end{lemma}

\begin{lemma} There exists a good selector $\vec{A}$. 
\begin{proof} Recall that $W_0=V(\G)\setminus V_\omega$. If $W_0$ is empty, then all selectors are good. Towards contradiction assume that all selectors are bad. In particular, $W_0$ is non-empty. Let $H$ be the set of increasing functions from $\kappa^+$ to $\N$. There are $2^{\kappa^+}$ many of them. For every $h\in H$ we will define a selector $\vec{A_h}$, in such a way that if $h\neq h'$, then $\vec{A_h}\nsim\vec{A_{h'}}$. This will give us $2^{\kappa^+}$ many bad selectors, which we saw is a contradiction.

We define $\vec{A_h}=\langle A_{h,\alpha},v_{h,\alpha},h(\alpha)\rangle$ simultaneously for all $h\in H$ by induction on $\alpha$. First let $A_{h,0}= W_0$, for all $h\in H$ and
$A_{h,\gamma}= \cap_{\alpha<\gamma} A_{h,\alpha}$ for  $\gamma$
limit. If $A_{h,\alpha}$ is given, by the assumption that all
selectors are bad, necessarily it must be non-empty. Choose
$v_{h,\alpha}$ in $A_{h,\alpha}$ and let
\[A_{h',\alpha+1}= \{v\in A_{h',\alpha}| \f{v,v_{h,\alpha}}=h'(\alpha)\},\]
for all $h'\in H$ such that $A_{h\rq{},\beta}=A_{h,\beta}$, for all $\beta\le\alpha$. In this way all selectors $\vec{A_{h\rq{}}}$ that agree with $\vec{A_{h}}$ up to $\alpha$ are assigned the \emph{same} element $v_{h,\alpha}$. The value $h\rq{}(\alpha)$ may be different and if $h'(\alpha)\neq h(\alpha)$, then, by the previous lemma, $\vec{A_{h'}}\nsim\vec{A_{h}}$. Since all functions differ from each other at a point, we get that all the $\vec{A_h}$ are non-equivalent. Contradiction.

Therefore, there exists a good selector.
\end{proof}
\end{lemma}

We will actually need slightly more than that. Assume that $H$ is the set of all increasing functions from $\kappa^+$ to $\N$ such that $h(\alpha)>m_{-1}$. There are still $2^{\kappa^+}$ many such functions and the argument of the previous lemma also proves that

\begin{corollary}\label{SelectorCorollary} There is a good selector $\vec{A}=\langle A_\alpha,v_\alpha,m_\alpha\rangle$ such that  for all $\alpha<\kappa^+$, $m_\alpha>m_{-1}$.
\end{corollary}

After this preliminary work on selectors, we are back to the goal of extending $f^{\G\rq{}}$ on pairs of the form $(b,v)$, for all $v\in W_0$ and for the unique $b\in V(\B)\setminus V(\A)$. 

\begin{definition}\label{FDefinition} Let $\vec{A}=\langle A_\alpha,v_\alpha,m_\alpha\rangle$ be  a good selector such that  $m_\alpha>m_{-1}$ for all $\alpha<\kappa^+$. Such a selector exists by corollary \ref{SelectorCorollary}.

Define \[\f{b,v_\alpha}=m_\alpha,\mbox{ for all $\alpha<\kappa^+$, and}\]
if $w\in A_\alpha\setminus A_{\alpha+1}$, then let \[\f{b,w}=\min\{\f{w,v_\alpha},\f{b,v_\alpha}=m_\alpha\}.\]
Since the selector is good, $\f{b,w}$ is defined for all $w\in W_0$.
\end{definition}
Definition \ref{FDefinition} together with definition \ref{VomegaDefinition} complete the definition of $f^{\G\rq{}}$ on pairs of the form $(b,v)$ for all $v\in V(\G)\setminus V(\A)$. Hence, $\G'$ is now a $\tau'$- structure. 

\begin{lemma}\label{PsiForGPrime} Under the above definition of $f$, any finitely generated/$\N$ substructure of $V(\G\rq{})=V(\G)\cup \{b\}$ satisfies sentence $\psi$  (cf. definition \ref{PsiDefinition}). In particular, any such substructure is isomorphic/$\N$ to a structure in $K(\N)$.
\begin{proof} Clauses (1) and (2) of $\psi$ are immediate. Clause (3) is the one that requires work. Since clause (3) involves only three (distinct) elements at a time, it is sufficient to consider substructures that are generated by three elements only. In particular, if all three elements are in $V(\G)$ the result is immediate by the assumptions on $\G$. So, assume $b$, the unique element of $V(\B)\setminus V(\A)$, is one of the three elements and let $a_1,a_2\in V(\G)$ be the other two. 

Before we proceed recall that for any $a\in V(\G)\setminus V(\A)$, $f(b,a)$ was defined by definition \ref{VomegaDefinition}, if $a\in V_\omega$, or by definition \ref{FDefinition}, if $a=v_\alpha$ or if $a\in A_{\alpha}\setminus A_{\alpha+1}$, some $\alpha<\kappa^+$. So, there are $3$ different cases we have to consider for each of $f^{\G'}(b,a_1)$ and  $f^{\G'}(b,a_2)$, resulting in a total of $9$ cases. Some of the cases will be consolidated together and the proof splits into corresponding cases. Our goal is to prove that the substructure/$\N$ generated by $b,a_1,a_2$ satisfies sentence $\psi$ and in particular, that if $f(b,a_1)\neq f(b,a_2)$, then $f(a_1,a_2)=\min\{f(b,a_1),f(b,a_2)\}$, otherwise $f(a_1,a_2)>f(b,a_1)=f(b,a_2)$.

If both $a_1,a_2$ are in $V_\omega$, then the substructure/$\N$ generated by $b,a_1,a_2$ satisfies sentence $\psi$ by claim \ref{VOmegaClaim}. Thus, we can assume that at least one of $a_1,a_2$ is not in $V_\omega$.

\begin{itemize}
  \item[Case I] There exists some $\alpha$ such that $a_1,a_2\in A_\alpha\setminus A_{\alpha+1}$.  

If one of $a_1,a_2$ is equal to $v_\alpha$, say $a_1$, then by definition \ref{FDefinition},  $\f{b,a_2}=\min\{\f{a_2,v_\alpha},\f{b,v_\alpha}\}=\min\{\f{a_2,a_1},\f{b,a_1}\}$, and the result is immediate. 

Hence, assume that both $a_1,a_2$ are not equal to $v_\alpha$ and let $x_1=\f{v_\alpha,a_1}, x_2=\f{v_\alpha,a_2}$. Since $a_1,a_2\notin A_{\alpha+1}$, $x_1,x_2\neq m_\alpha$. By definition \ref{FDefinition} again,  $\f{b,a_1}=\min\{x_1,m_\alpha\}$ and $ \f{b,a_2}=\min\{x_2,m_\alpha\}$.

\begin{itemize}
\item[Subcase I.1] $x_1\neq x_2$. 

Since $a_1,a_2,v_\alpha\in V(\G)$, the substructure/$\N$ generated by $a_1,a_2,v_\alpha$ satisfy $\psi$. In particular, since $x_1\neq x_2$, $\f{a_1,a_2}=\min\{x_1,x_2\}$. The reader can verify that in any linear ordering given three distinct elements, $x_1, x_2,m_\alpha$, exactly two out of the three quantities $\min\{x_1,m_\alpha\}$, $\min\{x_2,m_\alpha\}$ and $\min\{x_1,x_2\}$ are equal, while the third one is greater, which gives the result.

\item[Subcase I.2] $x_1=x_2<m_\alpha$. 

In this case $\min\{x_1,m_\alpha\}=x_1$ and $\f{b,a_1}$,  which is equal to $\min\{x_1,m_\alpha\}$, is also equal to $x_1$.  Similarly,  $\f{b,a_2}=\min\{x_2,m_\alpha\}=x_2$, i.e. $\f{b,a_1}=\f{b,a_2}$. In addition,  the substructure/$\N$ generated by $a_1,a_2,v_\alpha$ satisfies $\psi$ and since $x_1$ equals $x_2$, $f(a_1,a_2)$ must be greater than both $f(v_\alpha,a_1)=x_1$ and $f(v_\alpha,a_2)=x_2$. Overall, $\f{b,a_1}=\f{b,a_2}=x_1=x_2<\f{a_1,a_2}$. 

\item[Subcase I.3] $x_1=x_2>m_\alpha$. 

Then $\min\{x_1,m_\alpha\}=m_\alpha$ and $\f{b,a_1}$, which is equal to $\min\{x_1,m_\alpha\}$ is also equal to $m_\alpha$. Similarly, $\f{b,a_2}=\min\{x_2,m_\alpha\}=m_\alpha$. The same argument as  in Subcase I.2 for the substructure generated by $a_1,a_2,v_\alpha$ proves that $\f{a_1,a_2}$ is greater than both $f(v_\alpha,a_1)=x_1$ and $f(v_\alpha,a_2)=x_2$. Overall,$\f{b,a_1}=\f{b,a_2}=m_\alpha<x_1=x_2<\f{a_1,a_2}$.
\end{itemize}

\item[Case II] There exists some $\alpha$ such that $a_1\in A_\alpha\setminus A_{\alpha+1}$ and there exists some $\beta>\alpha$ such that $a_2=v_\beta$.

Let $\f{a_1,v_\alpha}=x$. Since $\beta>\alpha$, $A_\beta\subset A_{\alpha+1}$ and $v_\beta\in A_\beta$ implies that $v_\beta$ is also in $A_{\alpha+1}$. By definition of $A_{\alpha+1}$, $\f{v_\beta,v_\alpha}=m_\alpha$. Since $a_1,v_\alpha,v_\beta\in V(\G)$, the substructure/$\N$ generated by $a_1,v_\alpha,v_\beta$ satisfies $\psi$. So, $\f{a_1,v_\beta}=\min\{\f{a_1,v_\alpha},\f{v_\alpha,v_\beta}\}=\min\{x,m_\alpha\}$. By definition $\f{b,a_1}$ is equal to $\min\{x,m_\alpha\}$ and therefore, also equal to $\f{a_1,v_\beta}$. By definition also, $\f{b,v_\beta}=m_\beta$ and overall, $\f{b,v_\beta}=m_\beta>m_\alpha\ge \min\{x,m_\alpha\}=f(b,a_1)=f(a_1,v_\beta)$.

\item[Case III] $a_1,a_2,\alpha,\beta$ are as in case II, except that $\beta<\alpha$. 

We work similarly. In this case $A_\alpha\subset A_{\beta+1}$ and $a_1,v_\alpha\in A_\alpha$ implies that $a_1,v_\alpha$ are also in $A_{\beta+1}$. So,  $\f{v_\beta,a_1}=\f{v_\beta,v_\alpha}=m_\beta$ and by definition $\f{b,v_\beta}=m_\beta$. Considering the substructure/$\N$ generated by $v_\beta$, $v_\alpha$ and $a_1$, we conclude  that  $\f{v_\beta,a_1}=\f{v_\beta,v_\alpha}=m_\beta<\f{v_\alpha,a_1}$. Combine with the fact that $m_\beta<m_\alpha$ and we get that  $\f{b,v_\beta}=\f{v_\beta,a_1}=m_\beta<\min\{\f{v_\alpha,a_1},m_\alpha\}$. Since by definition $\min\{\f{v_\alpha,a_1},m_\alpha\}$ is equal to $\f{b,a_1}$, the result follows.

\item[Case IV] There exists $\alpha$ such that $a_1\in A_\alpha\setminus A_{\alpha+1}$ and $a_1\neq v_\alpha$,  and there exists some $\beta>\alpha$ such that $a_2\in A_\beta\setminus A_{\beta+1}$ and $a_2\neq v_\beta$.

First notice that since we can interchange $a_1$ and $a_2$, the same proof is good even for the case that $\beta<\alpha$.

Let $\f{a_1,v_\alpha}=x$ and $\f{a_2,v_\beta}=y$. Since $a_1\notin A_{\alpha+1}$, $x\neq m_\alpha$, and similarly, since $a_2\notin A_{\beta+1}$, $y\neq m_\beta$. The assumption $\beta>\alpha$ implies that $A_\beta\subset A_{\alpha+1}$ and $a_2,v_\beta\in A_\beta$ implies that $a_2,v_\beta$ are also in $A_{\alpha+1}$. So,  $\f{a_2,v_\alpha}=\f{v_\beta,v_\alpha}=m_\alpha$. Considering the substructure/$\N$ generated by $a_1,a_2,v_\alpha$ we conclude that $\f{a_1,a_2}=\min\{\f{a_1,v_\alpha},\f{v_\alpha,a_2}\}=\min\{x,m_\alpha\}$. By definition $\min\{x,m_\alpha\}$ is equal to $\f{b,a_1}$ and  $\f{b,a_2}$ equals $\min\{y,m_\beta\}$. If we can prove that $y>m_\alpha$, then $\f{b,a_2}=\min\{y,m_\beta\}>m_\alpha\ge\min\{x,m_\alpha\}=\f{a_1,a_2}=\f{b,a_1}$, which gives the result. In order to prove $y>m_\alpha$ consider the substructure/$\N$ generated by $a_2,v_\beta,v_\alpha$. Since $\f{a_2,v_\alpha}=\f{v_\beta,v_\alpha}=m_\alpha$, we get that $y=\f{a_2,v_\beta}>\f{a_2,v_\alpha}=\f{v_\beta,v_\alpha}=m_\alpha$. 

\item[Case V] There exists some $\alpha$ such that $a_1=v_\alpha$ and $a_2$ is in $V_\omega$. 

Since $a_1=v_\alpha\in W_0$, by the definitions of $V_\omega$ and $W_0$, $\f{b,a_2}=\f{a_2,v_\alpha}$. Also by definition $\f{b,v_\alpha}=m_\alpha$ and we assumed that for all $\alpha<\kappa^+$,     $m_\alpha>m_{-1}=\max\{\f{b,a}|a\in V_0\}$. By claim \ref{VZeroClaim}, $m_\alpha$ is also greater than $\max\{\f{b,a}|a\in V_\omega\}$. Therefore, $\f{b,a_2}=\f{a_2,v_\alpha}\le m_{-1}<m_\alpha=\f{b,v_\alpha}$.

\item[Case VI] There exists some $\alpha$ such that $a_1\in A_\alpha\setminus A_{\alpha+1}$ and $a_2$ is $V_\omega$.

Working similarly as in case V we can prove that  $\f{b,a_2}=\f{a_2,v_\alpha}<m_\alpha$. Moreover, by the definition of $W_0$, $\f{b,a_2}=\f{a_2,a_1}=\f{a_2,v_\alpha}$, and considering the substructure/$\N$ generated by $a_1,v_\alpha,a_2$, we conclude  $\f{a_2,a_1}=\f{a_2,v_\alpha}<\f{a_1,v_\alpha}$. Also, by definition,  $\f{b,a_1}=\min\{\f{a_1,v_\alpha},m_\alpha\}$. 

Combining all these,   $\f{b,a_1}=\min\{\f{a_1,v_\alpha},m_\alpha\}>\f{a_2,a_1}=\f{b,a_2}$.
  \end{itemize}
All these cases conclude the proof. 
\end{proof}
\end{lemma}

Lemma \ref{PsiForGPrime} proves that $\G\rq{}$ satisfies $(I)_{\N}$. By theorem \ref{towardsfullness}, there exists some $\G^*$ such that $\G^*$ is a model/$\N$, $|\G^*|=2^{(\kappa^+)}$ and $\G^*$ satisfies $\phi_\F$, the Scott sentence of $lim(K(\M))$ (cf. definition \ref{FraisseDefinition}). This concludes the proof of theorem \ref{SomeBig}. Since by theorem \ref{NotTooBig}, $\phi_\F$ can not have any models of size bigger than $2^{(\kappa^+)}$, $\phi_\F$ witnesses that $2^{(\kappa^+)}\in\ch$ which in return concludes the proof of the main theorem \ref{powersetofsuccessor}. 
\end{proof}
\end{theorem}

\begin{remark} \label{Phi0Remark}The only facts that we used in the proof of theorem \ref{powersetofsuccessor} about $\phi_0$ and the linear ordering $<$ are the following:
\begin{enumerate}
  \item The sentence $\phi_0$ is a complete sentence.
   \item In every model of $\phi_0$, $<$ does not have a maximum. 
  \item There is no model of $\phi_0$ of size $\kappa^{++}$.
  \item There is a model of $\phi_0$ with a $<$- cofinal sequence of length $\kappa^+$. (cf. Fact \ref{CofinalFact}). 
\end{enumerate}
\end{remark}

Replacing $\kappa^+$ by $\lambda$ and $\phi_0$ by $\phi$ in the proof of  theorem \ref{powersetofsuccessor}, we conclude:

\begin{theorem} Let $\lambda$ be an infinite cardinal and  $\phi$ be such that:
\begin{enumerate}
\item The sentence $\phi$ is a complete sentence.
  \item In every model of $\phi$, $<$ is a linear order without a maximum element.
  \item There is no model of $\phi$ of cardinality $\lambda^+$.
  \item There is a model of $\phi$ with a $<$- cofinal sequence of length $\lambda$.
\end{enumerate}
Then $2^\lambda$ is (homogeneously) characterizable.
\begin{proof} By the proof of theorem \ref{powersetofsuccessor} and remark \ref{Phi0Remark}.
\end{proof}
\end{theorem}

\begin{theorem}\label{powersetoflimit} Let $(\alephs{\alpha_n})_{n\in\omega}$ be an increasing sequence of cardinals in $\ch$ and $\alephs{\lambda}=\sup_n \alephs{\alpha_n}$. Then $2^{\alephs{\lambda}}$ is in \homch.
\begin{proof} By theorem \ref{powersetofsuccessor}, $2^{\alephs{\alpha_n+1}}$ is in \homch.

We know that \[2^{\alephs{\lambda}}=(2^{<\alephs{\lambda}})^{cf(\alephs{\lambda})},\]
and
\[2^{<\alephs{\lambda}}=\sup_n (2^{\alephs{\alpha_n+1}}).\]
By theorem \ref{limitthrm}, $2^{<\alephs{\lambda}}\in\ch$, and since $cf(\alephs{\lambda})=\omega$, we conclude that $2^{\alephs{\lambda}}\in\homch$ by theorem \ref{ltoomegahom}.
\end{proof}
\end{theorem}

Combining the last theorem with theorem \ref{powersetofsuccessor} we can conclude:

\begin{theorem}\label{mainhompowersetthrm} If $\alephs{\beta}\in\ch$, then
$2^{\alephs{\beta+\beta_1}}\in\homch$, for all
$0<\beta_1<\omegaone$.
\begin{proof}By induction on $\beta_1$. If it is a successor
ordinal, use theorem \ref{powersetofsuccessor}. If it is a
(countable) limit ordinal, then use \ref{powersetoflimit}.
\end{proof}
\end{theorem}

The only case that is not covered by the above theorem is when
$\alephs{\beta}$ is the head of a cluster (see definition
\ref{headofcluster}) and $\alephs{\beta}\neq \alephs{0}$ ($2^{\alephs{0}}$ is easily seen to be in $\homch$ by theorem \ref{ltoomegahom}). Combining with theorem
\ref{headpower} we get

\begin{theorem}\label{miscpowers}If $\alpha\le\beta$, $\alephs{\beta}\in\ch$, $\alpha_1<\omegaone$ and $0<\beta_1<\omegaone$, then
\[\alephs{\alpha+\alpha_1}^{\alephs{\beta+\beta_1}}\in\homch.\]
\begin{proof}By the previous theorem, $2^{\alephs{\beta+\beta_1}}\in\homch$. Since $\alpha\le\beta$, $\alephs{\alpha}^{\alephs{\beta+\beta_1}}=2^{\alephs{\beta+\beta_1}}$ and we conclude by theorem \ref{headpower}.
\end{proof}
\end{theorem}

Thus, depending on our model of ZFC, we get characterizability  of many cardinals that weren't
considered before, like $\alephalpha^{\alephs{\beta}}$ for
$\alpha,\beta<\omegaone$ etc.

\section{Open problems}
There are a few open problems of various difficulties:

\begin{enumerate}
  \item Is there any cardinal that is characterizable by a sentence in \lomegaone, but not characterizable by a Scott sentence? It is consistent that the answer is ``No''. Under the G.C.H. for instance, $\beths{\omegaone}$- the Hanf number for $\lomegaone$,- equals $\alephomegaone$ and all characterizable cardinals are strictly less than $\alephomegaone$. By \cite{HjorthsKnightPaper}, for all $\alpha<\omega_1$, $\alephalpha\in\ch$. These two results combined together imply that under the G.C.H. a cardinal $\kappa$ is characterizable iff $\kappa$ is completely characterizable iff $\kappa=\alephalpha$, for some $\alpha<\omega_1$.
 
So, we should ask if there is any model of ZFC in which the answer is positive. If there is no such model, then our job becomes a lot easier, since we do not have to worry about completeness every time. This is a big step forward. If there is such a model, then it will be very interesting to see one, but the author does not have any idea how such a model will look like. Either way, it seems to be a difficult question.
  \item Is there any cardinal in $\ch\setminus\homch$ other than $\aleph_0$? Can such a cardinal be a successor? If such a cardinal is a limit cardinal, does it necessarily have cofinality $\omega$? From \cite{CharacterizableCardinals} we know that it is consistent that all cardinals in $\ch\setminus\homch$ have cofinality $\omega$.

\begin{conjecture} For $\kappa$ a characterizable cardinal,
$\kappa$ is not homogeneously characterizable if and only if it has
cofinality $\omega$.
\end{conjecture}

  \item For $\kappa\in\ch$, does $2^\kappa\in\ch$? This would improve theorem \ref{mainhompowersetthrm} and would provide closure under powerset in all cases.
  If $\kappa^{\aleph_0}=\kappa$, then by theorem \ref{ltoomegahom},  $\kappa\in\homch$ and whence, $2^\kappa\in\homch$ by theorem \ref{baumgartner}. But does this hold for all $\kappa$? This of course depends on the model of ZFC too.
  \item If $\kappa,\lambda$ are both characterizable, does $\ltok\in\ch$? This would generalize the theorems we prove here. Under the G.C.H., the answer is trivially ``Yes", but is this a theorem in ZFC? By theorem \ref{headpower}, we have to consider only the case where $\lambda$ is the head of a cluster. See also \cite{LinearOrderings} for theorems along this line.
  \item Are there any closure properties for $\ch$ (and $\homch$) besides successor, countable unions, countable products, powerset and powers? 
 \end{enumerate}

\subsection*{Acknowledgement} The author would like to thank the anonymous referee for his/her persistence and dedication. 
The final version of this paper is a drastic improvement over the first version, which would have been impossible without the referee's comments.

\subsection*{In memory} This paper, as well as \cite{CharacterizableCardinals}, was written during the academic year 2006-2007, while visiting my thesis advisor, professor Greg Hjorth, at the University of Melbourne, Australia. Professor Hjorth died in 2011 and this paper is dedicated to his memory.

\bibliographystyle{amsplain}
\bibliography{Bibliography} 

\providecommand{\bysame}{\leavevmode\hbox to3em{\hrulefill}\thinspace}
\providecommand{\MR}{\relax\ifhmode\unskip\space\fi MR }
\providecommand{\MRhref}[2]{%
  \href{http://www.ams.org/mathscinet-getitem?mr=#1}{#2}
}
\providecommand{\href}[2]{#2}
\begin{thebibliography}{1}

\bibitem{BaumgartnersHanfNumber}
James~E. Baumgartner, \emph{{The Hanf number for complete
  $L_{\omega_1,\omega}$-sentences (without GCH).}}, J. Symb. Log. \textbf{39}
  (1974), 575--578 (English).

\bibitem{HanfsIncompactness}
W.~Hanf, \emph{{Incompactness in languages with infinitely long expressions.}},
  Fundam. Math. \textbf{53} (1964), 309--324 (English).

\bibitem{HjorthsKnightPaper}
Greg Hjorth, \emph{Knight's model, its automorphism group, and characterizing
  the uncountable cardinals}, J. Math. Log. \textbf{2} (2002), no.~1, 113--144.
  \MR{1900550 (2003c:03084)}

\bibitem{HodgesShorterModelTheory}
Wilfrid Hodges, \emph{{A shorter model theory.}}, {Cambridge: Cambridge
  University Press. x, 310 p.}, 1997 (English).

\bibitem{LaskowskiShelahsOnTheExistenceOfAtomicModels}
M.C. Laskowski and Saharon Shelah, \emph{{On the existence of atomic models.}},
  J. Symb. Log. \textbf{58} (1993), no.~4, 1189--1194 (English).

\bibitem{CharacterizableCardinals}
Ioannis Souldatos, \emph{Notes on cardinals that are characterizable by a
  complete (scott) sentence}, submitted, e-print can be found at
  http://arxiv.org/abs/1007.2426.

\bibitem{LinearOrderings}
\bysame, \emph{Linear orderings and powers of characterizable cardinals},
  Annals of Pure and Applied Logic \textbf{163} (2012), no.~3, 225 -- 237.

\end{thebibliography}

\end{document}